\documentstyle[12pt]{article}
\vspace*{-2cm}
\begin{document}
\begin{center}

\vskip1truecm

{\bf
SUPERALGEBRAS, THEIR QUANTUM
DEFORMATIONS\\ AND THE INDUCED REPRESENTATION METHOD} \footnote{~Work based on
preprint IC/96/130 and an invited talk given at the "{\it V--th National
Mathematics Conference}" (Hanoi, 17--20 September 1997).}\\[2mm]
({\it On the occasion of the 30--th anniversary of 
the Vietnam Mathematical Society})
\vskip1.5truecm
{\bf Nguyen Anh Ky}\\[4mm]
Institute of Physics\\
National Centre for Natural Science and Technology\\
P.O. Box 429, Bo Ho, Hanoi
10000,  Vietnam
\\[3mm] and\\[3mm]
International Centre for Theoretical
Physics, Trieste, Italy.
\end{center}
\vskip.5truecm
\centerline{ABSTRACT}
\baselineskip=20pt

Some introductory concepts and basic definitions of the Lie superalgebras
and their quantum deformations are exposed. Especially the induced
representation methods in both cases are described. Based on the Kac
representation theory we have succeeded in constructing representations
of several higher rank superalgebras. When representations of quantum
superalgebras are concerned, we have developed a method which can be
applied not only to the one--parametric quantum deformations but
also to the multi--parametric ones.  As an intermediate step the
Gel'fand--Zetlin basis description is extended to the case of
superalgebras and their quantum deformations. Our approach also
allows us to establish in consistent ways defining relations of quantum
(super)algebras. Some illustrations are given.

\vskip2truecm
\begin{center}
{\it Running title}: Superalgebras and quantum superalgebras.\\[1cm]
{\it Mathematical subject classification 1991}: 17A70, 81R50\\[1cm]
{\it Key words}: Superalgebras, quantum deformations, quantum
superalgebras, induced representations, Gel'fand--Zetlin basis.
\end{center}
\newpage
\begin{flushleft}
{\Large {\bf 1. Introduction}}
\end{flushleft}

The symmetry principles \cite{elot,wigner1,wigner2},
in particular, the supersymmetry idea \cite {susy1,susy2,susy3}, occupy
central places in
constructions of different physics theories. They are described by (super)
symmetry
groups or infinitesimally by corresponding (super) algebras. Especially,
superalgebras \cite{kac,kac3,kac4,sch} play an important role in
supersymmetry and
supergravity theories \cite {susy6,susy4,susy5}. They have various
applications in
quantum physics, superdynamical symmetry (superquantum mechanics), nuclear
physics, etc. Usually, as in the case of ordinary algebras, applications of
superalgebras lead us to finding explicit expressions for matrix elements
of their
generators. Therefore, constructing representations of superalgebras is an
actual
topic. Unfortunately, the purely mathematical problem is solved only partially.
Irrespective of the fact that all finite--dimensional irreducible
representations of the
basic classical Lie superalgebras are classified, the questions concerning
indecomposable finite--dimensional representations and constructions of
explicit
(indecomposable and irreducible) representations are less understood and
solved.
Especially for the orthosymplectic superalgebras it is not known how to
construct all
such indecomposable representations and modules. These (indecomposable and
irreducible) representations of the basic superalgebras and the structure
of the
corresponding modules were subjects of investigations of several authors, who
succeeded in constructing explicit representations and modules only for
lower rank
superalgebras \cite {rowe,pan,ritt}, while explicit representations of
larger
superalgebras were known very partially, besides some general expects \cite
{rowe,rowe2,vander}. Later, some movements forward were made in Refs.
\cite{k3,k1,k2,neli} where all finite--dimensional representations and a wide
class of
infinite--dimensional representations of several higher rank superalgebras
have already
been investigated in detail and constructed explicitly.\\

The quantum deformations \cite {drin,fa,jimbo,jimbo1,manin,wor},
originated from intensive investigations on the quantum scattering
problems and Yang--Baxter equations, represent an extension of the
symmetry concept. Since they became a subject of great interest, many
algebraic and geometric structures and some representations of quantum
(super-) groups and algebras have been obtained and understood (see in
this context, for example, Refs. \cite
{cele,chari,doeb,drin,fa,jimbo,jimbo1,kass,kulish,manin,man,wz,wor,wor2,yang}).
In particular, the quantum algebra $U_{q}[sl(2)]$ is very well studied
\cite{bieden,kure,macfa,pas,roche,skl}. As in the non--deformed case for
applications of quantum groups and
algebras we often need their explicit representations. However, although the
progress in this direction is remarkable the problem is still far from being
satisfactorily solved. Especially, representations of quantum superalgebras
\cite{chai,hoker,q-osc,k4,k6,k7,k8,k5,palev1,zhang} are presently under
development. Explicit
representations are known mainly for quantum superalgebras of lower
ranks and of particular
types like
$U_{q}[gl(n/1)]$, $U_{q}[osp(1/2)]$, etc., while for higher rank quantum
superalgebras of nonparticular type, only some general structures
\cite{zhang}, q--oscillator representations (see, for
example, Refs. \cite{hoker,q-osc}) and a class of representations of
$U_{q}[gl(m/n)]$ (Ref. \cite{palev1}) have been well investigated.
In general, representations, including the finite--dimensional
ones, of quantum superalgebras have not been explicitly constructed and
completely investigated (at neither generic $q$ nor $q$ being roots of
unity). Recently, in Ref. \cite{k4} we proposed an induced representation
method by which we can construct representations of higher rank quantum
superalgebras such as $U_{q}[gl(m/n)]$ for large $m,n$. \\[1mm]

 Here, in the framework of this paper, more precisely in the next
section, we shall make an introduction to the superalgebras and briefly
describe the induced representation method which is based on the
representation theory developed by Kac \cite{kac}. Then, in Sect. III, we
give a construction procedure for finite--dimensional representations of
the superalgebras $gl(m/n)$. The induced representation method allowed us
to construct explicitly all finite--dimensional representations and a wide
class of infinite--dimensional representations of several
higher rank superalgebras like $gl(2/2)$, $gl(3/2)$ and $osp(3/2)$. Sect.
IV is devoted to some introductory concepts of the quantum superalgebras.
Due to the method proposed in Ref. \cite{k4} and described in Sect. V we
succeeded for the first time in finding all finite--dimensional
representations, including the irreducible ones, of a higher rank quantum
superalgebra, namely $U_{q}[gl(2/2)]$ (see Refs. \cite{k4} and \cite{k5}).
It is clear that our method is applicable not only to other
one--parametric deformations but also to the multi--parametric ones
\cite{k6,k7,k8}.\\

Let us list some of the abbreviations and notations which will be used
throughout
the present paper:
\begin{tabbing}
\=12345678\=$V_{l}\otimes V_{r}$ -- \=tensor product between two linear spaces
$V_{l}$ and $V_{r}$\= or a tensor product \=\kill \>\>fidirmod(s) --
finite--dimensional
irreducible module(s), \\[2mm]
\>\>GZ basis -- Gel'fand--Zetlin basis,\\[2mm] \>\>lin. env.\{X\} -- linear
envelope of
X,\\[2mm] \>\>$q$ -- the deformation parameter,\\[2mm]
\>\>$[x]_{r}=(r^{x}-r^{-x})/(r-
r^{-1}), ~ r=r(q)$, ~ where $x$ is some number or operator,\\[2mm]
\>\>$[x]\equiv [x]_{q}$,\\[2mm]
\>\>$[E,F\}$ -- supercommutator between $E$ and $F$,\\[2mm]
\>\>$[E,F\}_{r}\equiv
EF\pm rFE$ -- r-deformed supercommutator between $E$ and $F$,\\[2mm]
\>\>$[m]$
-- a highest weight in a (GZ, for example,) basis $(m)$,\\[2mm]
\>\>$I_{k}^{q}$ -- the
maximal invariant subspace in $W^{q}([m])$, corresponding
to the class $k$,\\[2mm]
\>\>$W^{q}_{k}([m])=W^{q}([m])/I_{k}^{q}$ -- the class $k$ nontypical
module,\\[2mm]
\>\>$(m)^{\pm ij}$ -- a pattern obtained from $(m)$ by shifting $m_{ij}
\rightarrow
m_{ij} \pm 1$,\\[2mm] \>\> \end{tabbing}
Note that we must not confuse the quantum deformation $[x]\equiv [x]_{q}$
of $x$
with the highest weight (signature) $[m]$ in the GZ basis $(m)$ or with the
notation
$[ ~,~ ]$ for commutators. \\[9mm]
{\Large {\bf 2. Superalgebras and their representations}}\\

There exist several good references on Lie superalgebras and their
representations (see, for example, Refs. \cite{kac,kac3,kac4} and
\cite{sch}). Let us give here some introductory concepts and basic
definitions from the topic. A Lie
superalgebra
(from now on, only superalgebras) $A$, endowed with a Z-gradation, by
definition,
is a vector space which\\

1) is a direct sum of vector subspaces $A_{i}$, where $i \in {\bf Z}$:\\
$$A = \bigoplus_{i\in{\bf Z}}A_{i},
\eqno(2.1)$$ and\\

2) has a bilinear product (supercommutator) [ , \} such that

$$[x_{i},x_{j}\} := x_{i}x_{j} - (-1)^{ij}x_{j}x_{i} \in A_{i+j}, ~~ {\rm
for}~~ x_{i(j)} \in
A_{i(j)}, \eqno(2.2a) $$
$$[x,[y,z\}\} + [y,[z,x\}\} + [z,[x,y\}\} = 0. \eqno(2.2b)$$

One sees that the Lie superalgebra $A$ admits the following $Z_{2}$-graded
structure decomposition:

$$A = A_{\bar{0}} \bigoplus A_{\bar{1}}, \eqno(2.3a)$$ where
$$A_{\bar{0}(\bar{1})} = \bigoplus_{i=even(odd)}A_{i} \eqno(2.3b)$$ Here
$A_{\bar{0}}$, called the even subalgebra of $A$, is an ordinary Lie algebra,
while $A_{\bar{1}}$ is a subspace of the odd generators and represents an
$A_{\bar{0}}$--module as the supercommutator (2.2) defines in $A_{\bar{1}}$ a
homomorphism: $$A_{\bar{1}} \rightarrow A_{\bar{0}}\eqno(2.4)$$ We say the
above $Z$--gradation is consistent with the $Z_{2}$--one.\\[1mm]

Rewriting the decomposition (2.1a) in the form:

$$A = A_{-} \bigoplus A_{0} \bigoplus A_{+} \eqno(2.1a')$$ where
$$A_{-} = \bigoplus_{i<0}A_{i}, \hspace*{7mm} A_{+} = \bigoplus_{i>o}A_{i}
\eqno(2.1b')$$
we see that $A_{0}$,
referred to as a stability subalgebra, is either the even algebra
$A_{\bar{0}}$ or its
subalgebra and $A_{i}$'s are the adjoint representation spaces of $A$
restricted to
$A_{0}$ : $[A_{0},A_{i}]
\subseteq A_{i}$. The Cartan subalgebra is contained in $A_{0}$, while $A_{+}$
and $A_{-}$ are subspaces of the creation and annihilation generators,
respectively. One can construct a representation of $A$ induced from a
representation of the stability subalgebra $A_{0}$ by expressing the generators
from $A_{i}$ in a basis of the corresponding $Ad(A_{0})$--module.\\[1mm]

Let us denote by $V_{B}$ a module of a subalgebra $B$ of $A$. This $B$--module
$V_{B}$ can be extended to a $U(A)$--module, where $U(A)$ is the universal
enveloping algebra of $A$. An $A$--module $\tilde{V} := Ind_{B}^{A}V_{B}$
induced from the $B$--module is a $Z_{2}$--graded space obtained from $U(A)
\otimes V_{B}$ factorized by all the elements of the form $ab\otimes v - a
\otimes
b(v)$, $a\in A$, $b\in B$, $v\in V_{B}$ and endowed with the structure $a(u
\otimes
v) = au \otimes v$, $a\in A$, $u\in U(A)$, $v\in V_{B}$.\\

If $B = A_{0}\bigoplus A_{+}$ we can start from an $A_{0}$--module
$V_{0}(\Lambda)$, where
$\Lambda$ is a signature characterizing the corresponding representation of
$A_{0}$ in $V_{0}$. The latter becomes a $B$--module $V_{B}(\Lambda)$ by
setting
$$A_{+}V_{0} = 0 \eqno(2.5)$$
The induced module $\tilde{V}(\Lambda) = Ind_{B}^{A}V_{B}$, in general,
contains
a
(unique) maximal submodule $I(\Lambda)$.\\[7mm] {\bf Definition 2.1}:~ {\it
An
irreducible representation of a Lie superalgebra $A$ with the signature
$\Lambda$
is called the factor--module}
$$W(\Lambda) = \tilde{V}(\Lambda)/I(\Lambda,) \eqno(2.6)$$ {\it where $I$
is the
maximal submodule}.\\

Let

$$V = V_{\bar{0}} \bigoplus V_{\bar{1}} \eqno(2.7)$$ be a $Z_{2}$--graded
vector
space of the dimension $$dimV =(dimV_{\bar{0}},dimV_{\bar{1}}) =
(m,n)\eqno(2.8)$$ and $End(V)_{L}$ is a Lie
superalgebra of endomorphism mappings $End(V)$ endowed with the
multiplications (2.2).\\[7mm]
{\bf Definition 2.2}:~ {\it A linear representation of $A$ in $V$ is
defined as a homomorphism} $$\phi: A \rightarrow End(V)_{L} := gl(m/n).
\eqno(2.9)$$\\[1mm]

  We call $gl(m/n)$ a general linear Lie superalgebra which is a
super--analogue of the ordinary general linear Lie algebra $gl(m)$. Any
superalgebra is a subalgebra of $gl(m/n)$ and has the following matrix
representation:\\[1mm]
$$
\left[
\begin{array}{cc}
{\bf A} ~~~ & ~~~ {\bf C} \\[4mm]
{\bf D} ~~~ & ~~~ {\bf B}
\end{array}
\right]
\eqno(2.10)$$
where {\bf A, B, C} and {\bf D} are matrices of dimensions $m\times m,~
n\times n,~ m\times n$ and $n\times m$ , respectively. The even subalgebra
$A_{\bar{0}}$ is spanned by {\bf A}$\oplus ${\bf B} $\subseteq $ $gl(m)
\oplus gl(n)$, while {\bf C}
and {\bf D} are respectively the spaces of the positive-- and the negative
odd root
generators. For a basis of $gl(m/n)$ we can choose the Weyl matrices $e_{ij}$,
$$(e_{ij})_{kl}=\delta_{ik}\delta_{jl}, ~~~ i,j=1,2,...m+n, \eqno(2.11a)$$
satisfying the
supercommutation relations:\\[1mm] $$[e_{ij},e_{kl}\} = \delta_{jk}e_{il} -
(-1)^{[(i)+(j)][(k)+(l)]}\delta_{il}e_{jk}, ~~~ 1\leq i,j,k,l\leq m+n,
\eqno(2.11b)$$
where the gradation index $(i)$ is 0 for $1\leq i\leq m$ and 1 for $m+1\leq
i \leq
m+n$.\\[9mm]
{\Large {\bf 3. Representations of gl(m/n) in a gl(m) $\oplus$ gl(n) basis}}\\

Here we shall outline a construction procedure for a representation of the
superalgebra $gl(m/n)$ induced from a representation of the even subalgebra
$$A_{0}=gl(m/n)_{0}\equiv gl(m) \oplus gl(n)=A_{\bar{0}}\eqno(3.1)$$ in a
module
$V_{0}(\Lambda)$, where $\Lambda$ is some signature characterizing the
considered representation and being a highest weight in the case of a
finite--dimensional representation. The highest
weight $\Lambda$ represents an ordered set
$(\lambda_{1},\lambda_{2},...,\lambda_{m},\lambda_{m+1},...,\lambda_{m+n})$ of
the eigen--values $\lambda_{i}$ of the Cartan generators $e_{ii}$,
$i=1,2,..., m+n$,
on the so--called highest weight vector $M(\Lambda)$ which is defined as a
vector
from $V_{0}(\Lambda )$ and annihilated by the creation generators $e_{ij}$,
$1\leq
i<j\leq m$ or $m+1\leq i<j\leq m+n$, $$e_{ii}M=\lambda_{i}M, ~~~
\Lambda:=
(\lambda_{1},\lambda_{2},...,\lambda_{m},\lambda_{m+1},...,\lambda_{m+n}),
\eqno(3.2a)$$
$$e_{ij}M=0, ~~ {\rm for} ~~ 1\leq i<j\leq m ~~ {\rm or} ~~ m+1\leq i<j\leq m+n
\eqno(3.2b)$$

Identifying the subspaces $A_{\pm}$ as
$$A_{+} = \{e_{ij} ~\|~ m+n\geq j>m\geq i\geq 1\}, \eqno(3.3a)$$ $$A_{-} =
\{e_{ij} ~\|~ m+n\geq i>m\geq j\geq 1\} \eqno(3.3b)$$ we demand the condition
$$A_{+}V_{0}(\Lambda) = 0 \eqno(3.4)$$
which turns the module $V_{0}(\Lambda)$ in a $B$--module, where $$B =
gl(m)\oplus gl(n)\oplus A_{+}. \eqno(3.5)$$ The $gl(m/n)$--module $W(\Lambda)$
induced from the $gl(m) \oplus gl(n)$--module
$V_{0}(\Lambda )$ is the factor--space
$$W(\Lambda ) = (U \otimes V_{0}(\Lambda )/I(\Lambda )) \eqno(3.6)$$
where $U$ is the universal enveloping algebra of $gl(m/n)$, while
$I(\Lambda )$ is
the subspace
$$I(\Lambda ) ={\rm lin. env.}\left\{ ub\otimes v -u\otimes bv ~\|~
u\in U, b\in B\subset U, v\in V_{0}(\Lambda )\right\} \eqno(3.7)$$
The Poincar\'{e}--Birkhoff--Witt theorem states that $U$ is a linear span
of all the
elements
$$g = \prod_{e_{ij}\in A_{-}}(e_{ij})^{\theta_{ij}}b := a_{(-)}b,
\hspace*{6mm} b\in B,
\hspace*{4mm} \theta_{ij} = 0,1, \eqno(3.8)$$ where
$$a_{(-)}:=\prod_{e_{ij}\in A_{-}}(e_{ij})^{\theta_{ij}}\eqno(3.9)$$ are
ordered
sequences of the odd generators belonging to $A_{-}$. Now, considering
$g\otimes
v $ as an element of $W(\Lambda )$, from (3.9) we
have
$$g\otimes v = a_{(-)}b\otimes v
= a_{(-)}\otimes w, \hspace*{6mm} w=bv\in V_{0}(\Lambda )
\eqno(3.10)$$
Therefore
$$W(\Lambda) = {\rm lin. env.}\left\{ \prod_{e_{ij}\in
A_{-}}(e_{ij})^{\theta_{ij}}\otimes
v ~\|~ v\in V_{0}(\Lambda )\right\},
\hspace*{6mm}
{\it \theta_{ij} = 0,1} \eqno(3.11a)$$
or
$$W(\Lambda ) = T\otimes V_{0}(\Lambda) \eqno(3.11b)$$ where
$$ T = {\rm lin. env.}\left\{ \prod_{e_{ij}\in A_{-}}(e_{ij})^{\theta_{ij}}
~\|~
\theta_{ij}=0,1\right\}\subset U \eqno(3.12)$$
Since $T$, considered as an $Ad(A_{0})$--module, is $2^{mn}$--dimensional, the
module $W$ can be decomposed in a direct sum of a number ($2^{mn}$, at most)
of $A_{0}$--modules $V_{k}(\Lambda_{k})$ with highest weights $\Lambda_{k}$,
$0\leq k\leq (2^{mn}-1)$, i.e., $$W(\Lambda ) = \bigoplus_{0}^{2^{mn}-
1}V_{k}(\Lambda_{k}).\eqno(3.13)$$ where the notation
$$\Lambda_{0}\equiv \Lambda .\eqno(3.14)$$
is used. According to formulas (3.11), the
vectors $$\left |\theta_{ij};(m,n)\right >:=\prod_{e_{ij}\in
A_{-}}(e_{ij})^{\theta_{ij}}({\bf
m},{\bf n})= a_{(-)}({\bf m},{\bf n}) \eqno(3.15)$$
altogether span a basis of $W(\Lambda)$, where $({\bf m},{\bf n})$ is a
basis of
$V_{0}(\Lambda)$. Therefore, when $V_{0}$ is finite--dimensional, the
module $W$
and all other $gl(m/n)_{0}$--submodules $V_{k}$
are finite--dimensional, as well. For a basis of such a finite--dimensional
$gl(m/n)_{0}$--module we can choose the Gel'fand--Zetlin (GZ) tableaux also
called GZ (basis) vectors or patterns \cite{baird,barut,gz}:\\
$$({\bf m,n})_{k}\equiv \left(\left[
\begin{array}{c}
{\bf m}_{1m} ~ {\bf m}_{2m} ~~~ ...
~~~ {\bf m}_{m-1m} ~ {\bf m}_{mm}\\
{\bf m}_{1m-1} ~ {\bf m}_{2m-1} ~~ ...
~~ {\bf m}_{m-1m-1}\\
\vdots \\
{\bf m}_{12} ~ {\bf m}_{22}\\
{\bf m}_{11}
\end{array}
\right]
\otimes
\left[
\begin{array}{c}
{\bf n}_{1n} ~ {\bf n}_{2n} ~~~ ... ~~~
{\bf n}_{n-1n} ~ {\bf n}_{nn}\\
{\bf n}_{1n-1} ~ {\bf n}_{2n-1} ~~ ...
~~ {\bf n}_{n-1n-1}\\
\vdots\\
{\bf n}_{12} ~ {\bf n}_{22}\\
{\bf n}_{11}
\end{array}
\right]\right)_{k},\\[4mm]
\eqno(3.16)$$
where ${\bf m}_{ij}$ and ${\bf n}_{ij}$ are complex numbers, satisfying the
conditions:
$${\bf m}_{ij}-{\bf m}_{kl} \in {\bf Z},~~ {\bf m}_{ij} \geq {\bf m}_{ij-1}
\geq {\bf m}_{i-
1j}\eqno(3.17a)$$
and
$${\bf n}_{ij}-{\bf n}_{kl} \in {\bf Z},~~ {\bf n}_{ij} \geq {\bf n}_{ij-1}
\geq {\bf n}_{i-
1j},\eqno(3.17b)$$
as for $k=0$ we take
$$({\bf m}, {\bf n})_{0}\equiv ({\bf m}, {\bf n}). \eqno(3.18)$$ When there
does not
exist any threat of degenerations, the other subscripts $k$ are also
not necessary and therefore can be skipped. Thus, for every $V_{k}$ the highest
weight (signature) is characterized by the first row in (3.16)
$$\Lambda:=[\Lambda_{r},\Lambda_{l}]:= [{\bf m}_{1m}, {\bf m}_{2m}, ..., {\bf
m}_{mm}, {\bf n}_{1n}, {\bf
n}_{2n}, ..., {\bf n}_{nn}]:=[{\bf m},{\bf n}] \eqno(3.19)$$ combining the
highest
weights
$$\Lambda_{r}:=[{\bf m}]= [{\bf m}_{1m}, {\bf m}_{2m}, ..., {\bf m}_{mm}]
\eqno(3.20a)$$
and
$$\Lambda_{l}:=[{\bf n}]= [{\bf n}_{1n}, {\bf n}_{2n}, ..., {\bf n}_{nn}]
\eqno(3.20b)$$ of
$gl(m)$ and $gl(n)$, respectively.\\


As vectors from an $Ad(A_{0})$--module,
$a_{(-)}$ can be expressed in terms of a $gl(m) \oplus gl(n)$--GZ basis:
$$({\bf
m',n'})_{k}:=$$
\vspace*{4mm}
$$\left(\left[
\begin{array}{c}
{\bf m'}_{1m} ~ {\bf m'}_{2m} ~~~ ...
~~~ {\bf m'}_{m-1m} ~ {\bf m'}_{mm}\\
{\bf m'}_{1m-1}~ {\bf m'}_{2m-1} ~~ ...
~~ {\bf m'}_{m-1m-1}\\ \vdots\\
{\bf m'}_{12} ~ {\bf m'}_{22}\\
{\bf m'}_{11}
\end{array}
\right]
\otimes
\left[
\begin{array}{c}
{\bf n'}_{1n} ~ {\bf n'}_{2n} ~~~ ... ~~~ {\bf n'}_{n-1n} ~ {\bf n'}_{nn}\\
{\bf n'}_{1n-1} ~ {\bf n'}_{2n-1} ~~ ... ~~ {\bf n'}_{n-1n-1}\\
\vdots\\
{\bf n'}_{12} ~ {\bf n'}_{22}\\
{\bf n'}_{11}
\end{array}
\right]\right)_{k}
\eqno(3.21)$$
\vspace*{2mm}
Then, the basis (3.15) takes the form
$$\left |\theta_{ij};(m,n)\right >=\cal{N}.({\bf m'}, {\bf n'})\odot ({\bf
m},{\bf
n})\eqno(3.22)$$
where $\cal{N}$ is a norm.

The induced representations obtained, in general, are reducible in the
latter basis (3.22) referred to as an induced basis. In order to single
out all its irreducible subrepresentations we have to pass to another
basis, namely the reduced basis which is the union of all the GZ basis
vectors (3.16) for $k$ running from 0 to $2^{mn}-1$. The reduced basis
vectors (3.16) are connected with the induced basis ones (3.22) by the
Clebsch--Gordan decompositions written formally as follows
$$({\bf m}, {\bf n})_{k}=\sum_{(m,n),(m',n')}{\bf C}[(m,n)_{k}\|
(m',n');(m,n)]~({\bf
m'},{\bf n'})\odot ({\bf m},{\bf n})\eqno(3.23)$$ and vice versa
$$({\bf m'},{\bf n'})\odot ({\bf
m},{\bf n})=\sum_{(m,n)_{k}}{\bf C}^{-1}[(m',n');(m,n)\|(m,n)_{k}]~({\bf
m}, {\bf
n})_{k}\eqno(3.24)$$
where ${\bf C}$ and ${\bf C^{-1}}$ are short hands for the Clebsch--Gordan
coefficients and its invert expressions, respectively. The sums in (3.23)
and (3.24) spread over the Gel'fand--Zetlin ranges (3.17) for all possible
GZ patterns concerned. In Ref. \cite{k1} we proposed a modified GZ basis
description which can be extended later to the case of quantum superalgebras
\cite{k4,k6,k8,k5}.\\

   All matrix elements of $gl(m/n)$--generators in the induced basis or
in the reduced basis can be obtained by using formulas (2.11), (3.15),
(3.16), and (3.21)--(3.24). The main problem here is to find the
Clebsch--Gordan coefficients which are not always known explicitly,
especially for higher rank cases. For now, using the general method
described above, we can find all finite--dimensional representations of
the superalgebras $gl(2/2)$ and $gl(3/2)$, while the results for higher
rank $gl(m/n)$ are still partial.\\

   As an example we can consider the superalgebra $gl(2/1)$ generated by the
generators $e_{ij}$, $i,j=1,2,3$ satisfying (2.11) for $m=2,~ n=1$. Now the
space $T$ in (3.12) takes the following form $$T= {\rm lin. env}
\left\{(e_{31})^{\theta_{1}}(e_{32})^{\theta_{2}},
\theta_{i}=0,1.\right\}\eqno(3.25)$$
Then the module $W(\Lambda) $ (3.13) induced from a finite--dimensional
irreducible module (fidirmod) $V_{0}(\Lambda)$ of $gl(2/1)_0$ can be
decomposed
into four $gl(2/1)_0$--fidirmods $V_{k}$, $k=0,1,2,3$, $$W(\Lambda ) =
\bigoplus_{0}^{3}V_{k}(\Lambda_{k}).\eqno(3.26)$$ Now, the GZ basis
(3.16)--(3.18) for a finite--dimensional $gl(2/1)_0$--module $V_{k}$
represents a tensor
product $$\left[
\begin{array}{lcr}

\begin{array}{c}
m_{12}~~~m_{22}\\ m_{11}
\end{array}
;
\begin{array}{c}
m_{32}=m_{31}\\ m_{31}
\end{array}
\end{array}
\right]
\equiv
\left[
\begin{array}{lcr}

\begin{array}{c}
[m]_{2}\\ m_{11}
\end{array}
;
\begin{array}{c}
[m]_{1}\\ m_{31}
\end{array}
\end{array}
\right]
\equiv
(m)_{gl(2)}\otimes m_{31}\equiv (m)_{k}
\eqno(3.27a)$$
between the GZ basis $(m)_{gl(2)}$ of
$gl(2)$ and the $gl(1)$--factors $m_{31}$, where
$m_{ij}$ are complex numbers such that
$$m_{12}-m_{11},~ m_{11}-m_{22}\in
{\bf Z_{+}}\eqno(3.27b)$$ and
$$m_{32}=m_{31}.\eqno(3.27c)$$
Then $T$ as an $Ad(gl(2/1)_{0})$--module is spanned on the following basis
vectors
$$
1=
\left[
\begin{array}{lcr}

\begin{array}{c}
0~~~0\\ 0
\end{array}
;
\begin{array}{c}
0\\ 0
\end{array}
\end{array}
\right],\eqno(3.28a)$$
$$e_{3i}=
(-1)^i\left[
\begin{array}{lcr}

\begin{array}{c}
0~~~-1\\ i-2
\end{array}
;
\begin{array}{c}
1\\ 1
\end{array}
\end{array}
\right], ~~ i=1,2,\eqno(3.28b)$$
$$e_{31}e_{32}=
\left[
\begin{array}{lcr}

\begin{array}{c}
-1~~~-1\\ -1
\end{array}
;
\begin{array}{c}
2\\2
\end{array}
\end{array}
\right].\eqno(3.28c)$$

Using (3.11), (3.23), (3.24), (3.27) and (3.28) we can describe all the
basis vectors of the module $W$ and their transformations under the
actions of $gl(2/1)$--generators and then we can investigate the
irreducible and indecomposable structure of the module $W$. In such a way
all finite-dimensional irreducible representations of $gl(2/1)$ are found.
In order to make the present paper more compact we do not expose here
these results which represent classical limits of those of
$U_{q}[gl(2/1)]$ when $q\rightarrow 1$. The structure of the
$gl(2/1)$--module $W$ is similar to that of the module $W^q$ of
$U_{q}[gl(2/1)]$. The latter quantum superalgebra and its representations
will be considered (however, in a different approach) in section 5.\\

As far as the orthosymplectic Lie superalgebras $osp(m/n)$ are concerned, the
induced representation method is also applicable. However, this case has some
specific features which deserve to be mentioned. The orthosymplectic Lie
superalgebras $osp(m/n)$ which are a subclass of $gl(m/n)$ have various
applications in superfield theories \cite{sufield,susy6}, superquantum
mechanics \cite{suqm1,suqm2,suqm3,suqm4,suqm5}, nuclear physics
\cite{nphys,nphys2}, etc. Unfortunately, the mathematical
problem to determine the representations (or, say, only the finite--dimensional
representations) of the orthosymplectic Lie superalgebras is, at present,
far from
being solved. It is much less developed even in comparison to the other big
class of
the basis Lie superalgebras, namely $sl(m/n)$. Here, as an example, the
superalgebra $osp(3/2)$ is taken \cite{k2}.\\

In Ref. \cite{k2} we constructed explicitly all finite--dimensional
representations and a wide class of infinite--dimensional ones of $osp(3/2)$
induced from finite--representations of the stability algebra $A_{0}\equiv
so(3)\oplus gl(1)$ which is a subalgebra of the even subalgebra
$so(3)\oplus sp(2)$. The method depending on the representations of the
even algebras leads to an infinite--irreducible or indecomposable
$osp(3/2)$--module $\bar{W}(p,q)$ labeled by a
number pair ($p,q$). Any such module has, as mentioned above, a natural
induced basis, in which one easily writes transformations under the
actions of the generators. However, we need another basis, called reduced,
in order to easily single out and describe the invariant subspace
$\bar{W}_{inv}(p,q)$ of the module
$W(p,q)$ carrying infinite--irreducible or indecomposable
representations of $osp(3/2)$, the finite Kac module
$\bar{W}_{Kac}=\bar{W}/\bar{W}_{inv}$ (also carrying an irreducible or
indecomposable representation of $osp(3/2)$ and, finally, the irreducible
$osp(3/2)$ submodule $W_{Kac}(p,q)$ (which differs from the Kac module
only in the case of nontypical representations).\\[9mm]
{\Large {\bf 4. Quantum
superalgebras and their representations}}\\

As mentioned in the Introduction the quantum superalgebras and their
representations are subjects of intensive investigations in both physics and
mathematics. The quantum superalgebras as quantum deformations can be
introduced and defined in different ways
\cite{chai,hoker,drin,fa,q-osc,jimbo2,jimbo,k4,k6,k7,k8,k5,manin,man2,palev1,wor
,wor2,zhang}.

Here, we shall give some introductory concepts and basic definitions exposed
mostly in \cite{k4} where an induced representation method was proposed and
showed to be useful in constructing explicit representations of quantum
superalgebras \cite{k4,k6,k7,k8,k5}. Then, in the next section, for an
illustration of our
method we shall consider the quantum superalgebra $U{q}[gl(2/1)]$.\\[2mm]

Let $g$ be a rank $r$ (semi-) simple superalgebra, for example, $sl(m/n)$ or
$osp(m/n)$. The quantum superalgebra $U_{q}(g)$ as a quantum deformation
(q-deformation) of the universal enveloping algebra $U(g)$ of $g$, is
completely
defined by the Cartan-Chevalley canonical generators $h_{i}$, $e_{i}$ and
$f_{i}$,
$i=1,2,...,r$ which satisfy \cite{k4}\\[1mm]

1) {\bf the quantum Cartan--Kac supercommutation relations}
\begin{tabbing}
\=11111111111111111111111111111112\=$[h_{i},h_{j}]$ \= . =
0\=012345678901234567890123456789~\= \kill
\>\>$[h_{i},h_{j}]$\> = 0,\\[1mm]
\>\>$[h_{i},e_{j}]$\> = $a_{ij}e_{j}$,\\[1mm] \>\>$[h_{i},f_{j}]$\> =
$-a_{ij}f_{j}$,\\[1mm]
\>\>$[e_{i},f_{j}\}$\> = $\delta_{ij}[h_{i}]_{q_{i}^{2}}$,\>\>(4.1)
\end{tabbing}

2) {\bf the quantum Serre relations}
$$(ad_{q}{\cal E}_{i})^{1-\tilde{a}_{ij}}{\cal
E}_{j}=0,$$ $$(ad_{q}{\cal F}_{i})^{1-\tilde{a}_{ij}}{\cal
F}_{j}=0\eqno(4.2)$$ where
$(\tilde{a}_{ij})$ is a matrix obtained from the non-symmetric Cartan
matrix $(a_{ij})$
by replacing the strictly positive elements in rows with 0 on the diagonal
entry by $-
1$, while $ad_{q}$ is the q--deformed adjoint operator given by the formula
(4.8)\\[2mm] and\\[2mm]

3) {\bf the quantum extra--Serre relations}
(for $g$ being $sl(m/n)$ or $osp(m/n)$)
\cite{extra1,extra3,extra2}
$$\{[e_{m-1},e_{m}]_{q^{2}},[e_{m},e_{m+1}]_{q^{2}}\}=0,$$ $$\{[f_{m-
1},f_{m}]_{q^{2}},[f_{m},f_{m+1}]_{q^{2}}\}=0,\eqno(4.3)$$ being additional
constraints on the unique odd Chevalley generators $e_{m}$ and $f_{m}$. In the
above formulas we denoted $q_{i} = q^{d_{i}}$ where $d_{i}$ are rational
numbers
symmetrizing the Cartan matrix $d_{i}a_{ij}=d_{j}a_{ji}$, $1\leq i,j\leq
r$. For
example, in the case $g=sl(m/n)$ we have $$d_{i}=\left\{
\begin{array}{ll}
1 & ~~~~ {\rm if} ~~ 1\leq i\leq m,\\ -1 & ~~~~ {\rm if} ~~ m+1\leq i\leq
r=m+n-1.
\end{array}\right.
\eqno(4.4)$$

The above--defined quantum superalgebras form a subclass of a special class of
Hopf algebras called by Drinfel'd quasitriangular Hopf algebras \cite{drin}.
They are
endowed with a Hopf algebra structure given by the following additional maps:\\

1) {\bf coproduct} $\Delta$ : ~~$U$ $\rightarrow$ $U\otimes U$ \begin{tabbing}
\=11111111111111111111111111111\=$\Delta(h_{i})$\= $= h_{i}\otimes 1 +
1\otimes h_{i}$\=222222222222222\=123456789~\= \kill \>\>$\Delta(1)$\> $=
1\otimes
1$,\>\>\>\\[2mm] \>\>$\Delta(h_{i})$\> $= h_{i}\otimes 1 + 1\otimes
h_{i}$,\>\>\>\\[2mm] \>\>$\Delta(e_{i})$\> $= e_{i}\otimes q_{i}^{h_{i}} +
q_{i}^{-
h_{i}}\otimes e_{i}$,\>\>\>\\[2mm] \>\>$\Delta(f_{i})$\> $= f_{i}\otimes
q_{i}^{h_{i}} +
q_{i}^{-h_{i}}\otimes f_{i}$,\>\>\>(4.5) \end{tabbing}

2) {\bf antipode} $S$ : ~~~$U$ $\rightarrow$ $U$ \begin{tabbing}
\=11111111111111111111111123456\=$S(h_{i})$\= $=
-h_{i}$\=22222222222222222\=~12345678901201234\= \kill \>\>$S(1)$\> $=
1$,\>\>\>\\[2mm]
\>\>$S(h_{i})$\> $= -h_{i}$,\>\>\>\\[2mm] \>\>$S(e_{i})$\> $=
-q^{a_{ii}}_{i}e_{i}$,\>\>\>\\[2mm] \>\>$S(f_{i})$\> $=
-q^{-a_{ii}}_{i}f_{i}$ \>\>\>(4.6)
\end{tabbing}
\vspace*{2mm}
and\\

3) {\bf counit} $\varepsilon$ : ~$U$ $\rightarrow$ $C$ \begin{tabbing}
\=11111111111111111111111112345\=$S(h_{i})$\= $=
-h_{i}$1234567890123456789012345671234567~\= \kill \>\>$\varepsilon(1)$\> $=
1$,\>\\[2mm]
\>\>$\varepsilon(h_{i})$\>
$=\varepsilon(e_{i})=\varepsilon(f_{i})=0$,\>(4.7) \end{tabbing}
Then the quantum adjoint operator $ad_{q}$ has the following form
\cite{chai,roso}
$$ad_{q} = (\mu_{L} \otimes \mu_{R})(id \otimes S)\Delta\eqno(4.8)$$ with
$\mu_{L}$ (respectively, $\mu_{R}$) being the left (respectively, right)
multiplication:
$\mu_{L}(x)y = xy$ (respectively, $\mu_{R}(x)y = (-1)^{degx.degy} yx$).\\

A quantum superalgebra $U_{q}[gl(m/n)]$ is generated by the generators
$k_{i}^{\pm 1}\equiv q_{i}^{\pm E_{ii}}$, $e_{j}\equiv E_{j,j+1}$ and
$f_{j}\equiv
E_{j+1,j}$, $i=1,2,...,m+n$, $j=1,2,...,m+n-1$ such that the following
relations hold
\cite{k4}\\

1) {\bf the super--commutation relations}
\begin{tabbing}
\=1234567891234567
\=$k_{i}e_{j}k_{i}^{-1}$1\==1 \=$q_{i}^{(\delta_{ij}-
\delta_{i,j+1})}e_{j}$,~~~~
\=$k_{i}f_{j}k_{i}^{-1}$1\==12 \=$q_{i}^{(\delta_{ij+1}- \delta_{i,j})}f_{j}$,
12345678~\=\kill
\>\>$k_{i}k_{j}$\>=\>$k_{j}k_{i}$~, \>$k_{i}k_{i}^{-1}$\>=
\>$k_{i}^{-1}k_{i}$ =
1~,\\[1mm]
\>\>$k_{i}e_{j}k_{i}^{-1}$\>=\>$q_{i}^{(\delta_{ij} - \delta_{i,j+1})}e_{j}$~,
\>$k_{i}f_{j}k_{i}^{-1}$\>=\>$q_{i}^{(\delta_{ij+1} - \delta_{i,j})}f_{j}$~,
\\[1mm]
\>\>$[e_{i},f_{j}\}$\>=\> $\delta_{ij}[h_{i}]_{q^{2}_{i}}$, ~where
\>~~$q_{i}^{h_{i}}$\>=\>$k_{i}k_{i+1}^{-1}$,\>(4.9) \end{tabbing}

2) {\bf the Serre relations (4.2)} taking now the explicit forms
\begin{tabbing}
\=1234567891234567.\=$[e_{i},e_{j}]$\=~=~$[f_{i},f_{j}]$ \=~=~0,~~if
$|i-j|\neq 1$123\=2345678912345678~01234.\=\kill
\>\>$[e_{i},e_{j}]$\>~=~$[f_{i},f_{j}]$\>~=~0,~~if $|i-j|\neq 1$,\\[1mm]
\>\>~~$e_{m}^{2}$\>~=~~~$f_{m}^{2}$\>~=~0,\\[1mm]
\>\>$[e_{i},[e_{i},e_{j}]_{q^{\pm 2}}]_{q^{\mp 2}}$\>
\>~=~~$[f_{i},[f_{i},f_{j}]_{q^{\pm
2}}]_{q^{\mp 2}}$~=\>~0, ~~if $|i-j|=1$\>(4.10)\\[1mm] \end{tabbing}
\vspace*{2mm}
and\\

3) {\bf the extra--Serre relations (4.3)}
$$\{[e_{m-1},e_{m}]_{q^{2}},[e_{m},e_{m+1}]_{q^{2}}\}~= 0,$$ $$\{[f_{m-
1},f_{m}]_{q^{2}},[f_{m},f_{m+1}]_{q^{2}}\}~= 0.\eqno(4.3)$$
Here, besides $d_{i}$, $1\leq i\leq r=m+n-1$ given in (4.4), we introduced
$$d_{m+n}=-1.\eqno(4.11)$$
The Hopf structure on $k_{i}$ looks like \begin{tabbing}
\=1234567891234567891234567891234\=$\Delta(k_{i})$~\= =~ \=$k_{i}\otimes
k_{i}$123456789123456789123456789\=\kill
\>\>$\Delta(k_{i})$\> = \>$k_{i}\otimes k_{i}$,\\[1mm] \>\>$S(k_{i})$\> =
\>$k_{i}^{-
1}$,\\[1mm] \>\>$\varepsilon(k_{i})$\> = \>1.\>(4.12) \end{tabbing}
The generators $E_{ii}$, $E_{i,i+1}$ and $E_{i+1,i}$ together with the
generators
defined in the following way \begin{tabbing}
\=1234567890123\=$E_{i,j+1}$\=:= $[E_{ik}E_{kj}]_{q^{-2}}$x\= $\equiv
E_{ik}E_{kj}~-~q^{-2}E_{kj}E_{ik}$,
\=$~~i<k<j$,\=123456789\=\kill
\>\>$E_{ij}$\>:= ~$[E_{ik}E_{kj}]_{q^{-2}}$\> $\equiv ~E_{ik}E_{kj}~-~q^{-
2}E_{kj}E_{ik}$,\>$~~i<k<j$, \\[1mm]
\>\>$E_{ji}$\>:= ~$[E_{jk}E_{ki}]_{q^{2}}$\> $\equiv ~E_{jk}E_{ki}~-
~q^{2}E_{ki}E_{jk}$$,~$\>$~~i<k<j$,\>\> (4.13) \end{tabbing}
play an analogous role as the Weyl generators $e_{ij}$, $$(e_{ij})_{kl} =
\delta_{ik}\delta_{jl}, ~ i,j=1,2,...,m+n\eqno(2.11a)$$ of the superalgebra
$gl(m/n)$
whose universal enveloping algebra $U[gl(m/n)]$ represents a classical limit of
$U_{q}[gl(m/n)]$ when $q\rightarrow 1$.\\

The quantum algebra $U_{q}[gl(m/n)_{0}]\cong U_{q}[gl(m)\oplus gl(n)]$
generated
by $k_{i}$, $e_{j}$ and $f_{j}$, $i=1,2,...,m+n$, $m\neq j=1,2,...,m+n-1$,
$$U_{q}[gl(m/n)_{0}] ~=~{\rm lin. env.}\{E_{ij}\|~ 1\leq i,j\leq m ~~and~~
m+1\leq
i,j\leq m+n\}\eqno(4.14)$$ is an even subalgebra of $U_{q}[gl(m/n)]$. Note that
$U_{q}[gl(m/n)_{0}]$ is included in the largest even subalgebra
$U_{q}[gl(m/n)]_{0}$ containing all elements of $U_{q}[gl(m/n)]$ with even
powers
of the odd generators.\\

In Ref. \cite{k4} we describe the construction method for induced
representations of
$U_{q}[gl(m/n)]$ and give in detail a procedure of how to construct all
finite--dimensional representations of $U_{q}[gl(2/2)$ (see Refs. \cite{k4}
and \cite{k5}).
Let us briefly explain why and how we can use the proposed induced
representation method, which can be applied not only to one--parametric
deformations but also to multi--parametric ones \cite{k6,k7,k8}.\\


Indeed, our method is based on the fact \cite{lus,roso2} that a
finite--dimensional
representation of a Lie algebra $g$ can be
deformed to a finite--dimensional representation of its quantum analogue
$U_{q}(g)$. In particular, finite--dimensional representations of
$U_{q}[gl(m)\oplus
gl(n)]$ are quantum deformations of those of $gl(m)\oplus gl(n)$. This
means that a
finite--dimensional irreducible representation of $U_{q}[gl(m)\oplus gl(n)]$
is again
highest weight. On the other hand, as we can see from (4.3), (4.9)--(4.11),
(4.13)
and (4.14), $U_{q}[gl(m)\oplus gl(n)]$ is
the stability subalgebra of $U_{q}[gl(m/n)]$. Therefore, we can construct
representations of $U_{q}[gl(m/n)]$ induced from finite--dimensional
irreducible
representations of $U_{q}[gl(m)\oplus gl(n)]$. Let $V_{0}^{q}(\Lambda)$ be a
$U_{q}[gl(m)\oplus gl(n)]$--fidirmod characterized by some highest weight
$\Lambda$. The module $V_{0}^{q}(\Lambda)$ represents a tensor product
between a $U_{q}[gl(m)]$--fidirmod $V_{0,m}^{q}(\Lambda_{m})$ of a highest
weight $\Lambda_{m}$ and a $U_{q}[gl(n)]$--fidirmod $V_{0,n}^{q}(\Lambda_{n})$
of a highest weight $\Lambda_{n}$
$$V_{0}^{q}(\Lambda)=V_{0,m}^{q}(\Lambda_{m})\otimes
V_{0,n}^{q}(\Lambda_{n}) \eqno(4.15)$$
where $(\Lambda_{m})$ and $(\Lambda_{n})$ are defined respectively as the
left and right components of $\Lambda$
$$\Lambda =[\Lambda_{m}, \Lambda_{n}].\eqno(4.16)$$ For a basis of each of
$V_{0,m}^{q}(\Lambda_{m})$ and $V_{0,n}^{q}(\Lambda_{n})$, i.e., of
$V_{0}^{q}(\Lambda)$ we can choose the Gel'fand-Zetlin (GZ) tableaux
\cite{barut,baird,gz}, since the latter are invariant under the quantum
deformations \cite{chere,jimbo2,lus,roso2,tolst,ueno}. Therefore,
the highest weight $\Lambda$ is
described again by the first rows of the GZ tableaux called from now on as
the GZ (basis) vectors.\\

Demanding $$E_{m,m+1}V_{0}^{q}(\Lambda)\equiv
e_{m}V_{0}^{q}(\Lambda)=0\eqno(4.17)$$ i.e.
$$U_{q}(A_{+})V_{0}^{q}(\Lambda)=0 \eqno(4.18)$$ we turn $V_{0}^{q}(\Lambda)$
into a $U_{q}(B)$-module, where $$A_{+} = \{E_{ij}\|~ 1\leq i\leq m < j\leq
m+n\},
\eqno(4.19)$$ $$B = A_{+}\oplus gl(m)\oplus gl(n). \eqno(4.20)$$ The
$U_{q}[gl(m/n)]$ -module $W^{q}$ induced from the
$U_{q}[gl(m)\oplus gl(n)]$-module $V_{0}^{q}(\Lambda)$ is the factor-space
$$W^{q}=W^{q}(\Lambda) =
[U_{q} \otimes
V_{0}^{q}(\Lambda)]/I^{q}(\Lambda)
\eqno(4.21)$$
where $U_{q}\equiv U_{q}[gl(m/n)]$, while $I^{q}(\Lambda)$ is the subspace
$$I^{q}(\Lambda) = {\rm lin. env.} \left\{ ub\otimes v -u\otimes bv \|
~ u\in U_{q},~ b\in U_{q}(B)\subset U_{q},~ v\in V_{0}^{q}(\Lambda)\right\}.
\eqno(4.22)$$

Any vector $w$ from the module $W^{q}$ has the form $$w=u\otimes v,~~~ u\in
U_{q},~~ v\in V_{0}^{q}\eqno(4.23)$$ Then $W^{q}$ is a
$U_{q}[gl(m/n)]$-module in
the sense $$gw\equiv g(u\otimes v)=gu\otimes v\in W^{q}\eqno(4.24)$$ for $g$,
$u\in U_{q}$, $w\in W^{q}$ and $v\in V_{0}^{q}$.\\

As we can see from (4.17) the modules $W^{q}(\Lambda)$ and
$V_{0}^{q}(\Lambda)$ have one and the same highest weight vector. Therefore,
they are characterized by one and
the same highest weight $\Lambda$.\\[9mm]
\noindent
{\Large {\bf 5. Induced representations
of U$_{q}$[gl(2/1)]}}\\

Although the general expressions of the finite--dimensional representations
of the quantum superalgebra $U_{q}[gl(2/1)]$ can be found from
\cite{palev2,zhang} for $n=2$, the irreducible representations,
however, have not yet been considered in detail. Now, as an illustration
of the method described above \cite{k4}, we investigate and construct
explicitly all irreducible (i.e., typical and nontypical) finite--dimensional
representations of
$U_{q}[gl(2/1)]$. Here, we assume that the quantum deformation parameter $q$ is
generic. It means that there does not exist any positive integer $N\in {\bf
Z^+}$ such that $q^{N}=1$. We can construct directly and explicitly
representations of the quantum superalgebra $U_{q}[gl(2/1)]$ induced from
some (usually, irreducible) finite--dimensional representations of the
even subalgebra $U_{q}[gl(2)\oplus gl(1)]$. Since the latter is a
stability subalgebra of $U_{q}[gl(2/1)]$ we expect that
the constructed induced representations of $U_{q}[gl(2/1)]$ are decomposed
into finite--dimensional irreducible representations of $U_{q}[gl(2)\oplus
gl(1)]$. For this purpose we shall introduce a $U_{q}[gl(2/1)]$--basis
(i.e., a basis within a $U_{q}[gl(2/1)]$--module or shortly, a basis of
$U_{q}[gl(2/1)]$) convenient for us in investigating the module structure.
This basis (see (5.29)) can be expressed in terms of some basis of the
even subalgebra $U_{q}[gl(2)\oplus gl(1)]$ which in turn
represents a (tensor) product between a $U_{q}[gl(2)]$--basis and a
$gl(1)$--factor. It will be shown that the finite--dimensional
representations of $U_{q}[gl(2)]$, i.e., of
$U_{q}[gl(2)\oplus gl(1)]$
can be realized in the Gel'fand--Zetlin (GZ) basis. The finite--dimensional
representations of $U_{q}[gl(2/1)]$ constructed are irreducible and can be
decomposed into finite--dimensional irreducible representations of the
subalgebra $U_{q}[gl(2)\oplus gl(1)]$.\\

The quantum superalgebra $U_{q}[gl(2/1)]$ is completely defined through the
Cartan--Chevalley generators $E_{12}$, $E_{21}$, $E_{23}$, $E_{32}$, and
$E_{ii}$, $i=1,2,3$, satisfying the relations (4.9) and (4.10) which now read
\begin{tabbing}
\=12345678912345678912345678\=$[E_{ii},E_{jj}]$123\= =12\= 0123
\=$[E_{ii},E_{j,j+1}]$\= =
\=$(\delta_{ij}-\delta_{i,j+1})E_{j,j+1}$\=1234567\kill

~~~~~~1) {\bf the super--commutation relations} ($1\leq i,i+1,j,j+1\leq
3$):\\[2mm]
\>\>$[E_{ii},E_{jj}]$\> = \>0,\>\>\>\>(5.1a)\\[1mm]
\>\>$[E_{ii},E_{j,j+1}]$\> =
\>$(\delta_{ij}-\delta_{i,j+1})E_{j,j+1}$, \>\>\>\>(5.1b)\\[1mm]
\>\>$[E_{ii},E_{j+1,j}]$\> =
\>$(\delta_{i,j+1}-\delta_{ij})E_{j+1,j}$, \>\>\>\>(5.1c)\\[1mm]
\>\>$[E_{12},E_{32}]$\> = \> $[E_{21},E_{23}]~=~0$, \>\>\>\>(5.1d)\\[1mm]
\>\>$[E_{12},E_{21}]$\> = \>$[h_{1}]$,\>\>\>\>(5.1e)\\[1mm]
\>\>$\{E_{23},E_{32}\}$\>
= \>$[h_{2}]$,\>\>\>\>(5.1f)\\[1mm] \>\>$h_{i}$\>=\>$(E_{ii}-{d_{i+1}\over
d_{i}}E_{i+1,i+1}),$\>\>\>\>(5.1g)\\[4mm] with $d_{1}=d_{2}=-d_{3}=1$,
\end{tabbing}
\vspace*{6mm}

~~~~~~2) {\bf the Serre--relations}:
$$E_{23}^{2}=E_{32}^{2}=0,$$
$$[E_{12},E_{13}]_{q}=0,$$
$$[E_{21},E_{31}]_{q}=0, \eqno(5.2)$$
respectively, where
$$E_{13}:=[E_{12},E_{23}]_{q^{-1}}\eqno(5.3a)$$ and
$$E_{31}:= -[E_{21},E_{32}]_{q^{-1}}. \eqno(5.3b)$$ are defined as new odd
generators which have vanishing squares. Now the extra--Serre relations are not
necessary, unlike higher rank cases \cite{extra1,extra3,k4,k5,extra2}.\\

 As mentioned earlier, these generators $E_{ij}$, $i,j= 1,2,3$, are quantum
deformation analogues (q--analogues) of the Weyl generators $e_{ij}$
$$(e_{ij})_{kl}=\delta_{ik}\delta_{jl},~~i,j,k,l=1,2,3,\eqno(5.4)$$ of the
classical (i.e., non--deformed) superalgebra $gl(2/1)$ whose universal
enveloping algebra $U[gl(2/1)]$ is a classical limit of $U_{q}[gl(2/1)]$
when $q\rightarrow 1$.\\

>From the relations (5.1)--(5.3) we see that every of the odd spaces
$A_{\pm}$ $$A_{+}= {\rm lin. env.}\{E_{13},E_{23}\}, \eqno(5.5)$$
$$A_{-}= {\rm lin. env.}\{E_{31},E_{32}\},\eqno(5.6)$$ as always, is a
representation
space of the even subalgebra $U_{q}[gl(2/1)_{0}]\equiv U_{q}[gl(2)\oplus
gl(1)]$
which, generated by the generators $E_{12}$, $E_{21}$, and $E_{ii}$,
$i=1,2,3$, is
a stability subalgebra of $U_{q}[gl(2/1)]$. Therefore, we
can construct representations of $U_{q}[gl(2/1)]$ induced from some (usually,
irreducible) representations of $U_{q}[gl(2/1)_{0}]$ which are realized in some
representation spaces (modules) $V^{q}_{0}$ being tensor products of
$U_{q}[gl(2)]$--modules
$V^{q}_{0, gl_{2}}$ and $gl(1)$--modules (factors) $V^{q}_{0,gl_{1}}$.
Following Ref. \cite{k4} we demand
$$E_{23}V_{0}^{q}=0\eqno(5.7)$$
that is
$$U_{q}(A_{+})V_{0}^{q}=0.\eqno(5.8)$$
In such a way we turn the $U_{q}[gl(2/1)_{0}]$--module $V^{q}_{0}$ into a
$U_{q}(B)$--module where
$$B=A_{+}\oplus gl(2)\oplus gl(1).\eqno(5.9)$$ The $U_{q}[gl(2/1)]$--module
$W^{q}$ induced from the $U_{q}[gl(2/1)_{0}]$--module $V^{q}_{0}$ is the
factor--space $$W^{q}=[U_{q}\otimes V_{0}^{q}]/I^{q}\eqno(5.10)$$ where
$$U_{q}\equiv U_{q}[gl(2/1)],\eqno(5.11)$$ while $I^{q}$ is the subspace
$$I^{q}={\rm lin. env.}\{ub\otimes v-u\otimes bv\| u\in U_{q}, b\in
U_{q}(B)\subset
U_{q}, v\in V_{0}^{q}\}.\eqno(5.12)$$

Any vector $w$ from the module $W^{q}$ is represented as $$w=u\otimes v,~~~~
u\in U_{q},~~~~ v\in V_{0}^{q}.\eqno(5.13)$$ Then $W^{q}$ is a
$U_{q}[gl(2/1)]$--module in the sense $$gw\equiv g(u\otimes v)=gu\otimes
v\in W^{q}\eqno(5.14)$$
for $g,~u\in U_{q}$, $w\in W^{q}$ and $v\in V_{0}^{q}$.\\

Using the commutation relations (5.1)--(5.2) and the definitions (5.3) we
can prove
the following analogue of the Poincar\'e--Birkhoff--Witt
theorem\\[4mm]
{\bf Lemma 5.1}: {\it The quantum deformation $U_{q} := U_{q}[gl(2/1)]$ is
spanned on all possible linear combinations of the elements} $$g =
(E_{23})^{\eta_{1}}(E_{13})^{\eta_{2}}(E_{31})^{\theta_{1}}
(E_{32})^{\theta_{2}}g_{0},\eqno(5.15)$$ {\it where $\eta_{i}$,
$\theta_{i}=0,1$ and
$g_{0}\in U_{q}[gl(2/1)_{0}]\equiv U_{q}[gl(2)\oplus gl(1)]$}.\\

Taking into account (5.10)--(5.12) and (5.15) we arrive at the following
assertion\\[4mm]
{\bf Lemma 5.2}: {\it The induced $U_{q}[gl(2/1)]$--module $W^{q}$ is the
linear span}
$$W^{q}={\rm
lin. env.}\{(E_{31})^{\theta_{1}}(E_{32})^{\theta_{2}}\otimes v\|v\in
V_{0}^{q},~~\theta_{1}, ~\theta_{2}=0,1\},\eqno(5.16)$$ {\it and,
consequently, all
the vectors of the form} $$\left |\theta_{1}, \theta_{2}; (m)\right > :=
(E_{31})^{\theta_{1}}(E_{32})^{\theta_{2}}\otimes (m), ~~ \theta_{1},~
\theta_{2}=0,1, \eqno(5.17)$$
{\it constitute a basis in $W^{q}$, where $(m)$ is a (GZ, for example,)
basis in
$V_{0}^{q}$}.\\

  Therefore, if $V_{0}^{q}$ is a finite--dimensional
$U_{q}[gl(2/1)_{0}]$--module, the
$U_{q}[gl(2/1)]$--module $W^{q}$ is finite--dimensional, as well. Moreover,
$W^{q}$ is a highest weight module due to the condition (5.7) imposed on
$V_{0}^{q}$ which, as a finite--dimensional $U_{q}[gl(2/1)_{0}]$--module,
is always highest weight. Then, based on the latest {\bf Lema 5.2} and the
fact that
$U_{q}[gl(2/1)_{0}]$ is a stability subalgebra of $U_{q}[gl(2/1)]$ we
conclude that
$W^{q}$ can be decomposed in a number of $U_{q}[gl(2/1)_{0}]$--modules which
in turn are finite--dimensional and, therefore, highest weight. Such a
highest weight
module is characterized by a signature (referred otherwise to as a highest
weight), being an ordered set of the eigen--values of the Cartan
generators on the so--called
highest weight vector defined as a vector annihilated by the creation
generators (see (3.4) and (3.5)). Thus, the condition (5.7) means that
$W^{q}$ and $V_{0}^{q}$
have one and the same highest weight vector, i.e., one and
the same highest weight. Let $V_{0}^{q}$ be a finite-dimensional irreducible
module (fidirmod) of $U_{q}[gl(2/1)_{0}]$.\\[4mm]
{\bf Lemma 5.3}: {\it The
$U_{q}[gl(2/1)]$--module $W^{q}$ is decomposed into (4 or less)
finite--dimensional irreducible
modules $V_{k}^{q}$ of the even subalgebra $U_{q}[gl(2/1)_{0}]$}
$$W^{q}([m])=\bigoplus _{0\leq k\leq 3}V_{k}^{q}([m]_{k}), \eqno(5.18)$$
{\it where $[m]$ and $[m]_{k}$ are some signatures (highest--weights)
characterizing the module $W^{q}\equiv W^{q}([m])$ and the modules
$V_{k}^{q}\equiv V_{k}^{q}([m]_{k})$, respectively}.\\
The proof of this lemma follows from thr very construction of $W^q$ and
the same argument used for deriving (3.13).\\

  Each of the fidirmods $V_{k}^{q}$,
$0\leq k\leq 3$, is spanned on a basis, say $(m)_{k}$, which can be taken
as a tensor product between a GZ basis of $U_{q}[gl(2)]$ and $gl(1)$--factors.
In this case, we also call $(m)_{k}$ as a GZ basis. It is clear that
$$(m)_{0}\equiv (m)\eqno(5.19a)$$
and
$$[m]_{0}\equiv [m]\eqno(5.19b)$$
in our notations. We refer to
the basis (5.17) as the
induced $U_{q}[gl(2/1)]$--basis (or simply, the induced basis) in order to
distinguish
it from the reduced one introduced later in the next subsection.
\begin{flushleft}
{\bf 5.a. Finite--dimensional representations of U$_{q}$[gl(2/1)]}\\
\end{flushleft}

We can show that finite--dimensional representations of
$U_{q}[gl(2/1)_{0}]$ can be
realized in some spaces (modules) $V_{k}^{q}$
spanned by the (tensor) products
$$\left[
\begin{array}{lcr}

\begin{array}{c}
m_{12}~~~m_{22}\\ m_{11}
\end{array}
;
\begin{array}{c}
m_{32}=m_{31}\\ m_{31}
\end{array}
\end{array}
\right]
\equiv
\left[
\begin{array}{lcr}

\begin{array}{c}
[m]_{2}\\ m_{11}
\end{array}
;
\begin{array}{c}
[m]_{1}\\ m_{31}
\end{array}
\end{array}
\right]
\equiv
(m)_{gl(2)}\otimes m_{31}\equiv (m)_{k}
\eqno(5.20a)$$
between the (GZ) basis vectors $(m)_{gl(2)}$ of $U_{q}[gl(2)]$ and the
$gl(1)$--factors $m_{31}$, where $m_{ij}$ are complex numbers such that
$$m_{12}-m_{11},~ m_{11}-m_{22}\in
{\bf Z_{+}}\eqno(5.20b)$$ and
$$m_{32}=m_{31}.\eqno(5.20c)$$
Indeed, any finite--dimensional representation of (not only) $U_{q}[gl(2)]$
is always
highest weight and if the generators $E_{ij}$, $i,j=1,2$ and $E_{33}$ are
defined on
(5.20) as follows \begin{eqnarray*}
~~~~~~~~~~~~~~~~~~~~~~~~~~~ E_{11}(m)_{k}& = &(l_{11}+1)(m)_{k},\\
E_{22}(m)_{k}& = &(l_{12}+l_{22}-l_{11}+2)(m)_{k},\\ E_{12}(m)_{k}& =
&\left([l_{12}-l_{11}][l_{11}-l_{22}]\right)^{1/2}(m)_{k}^{+11},\\
E_{21}(m)_{k}& =
&\left([l_{12}-l_{11}+1][l_{11}-l_{22}-1]\right)^{1/2}(m)_{k}^{-11},\\
E_{33}(m)_{k}& =
& (l_{31}+1)(m)_{k},
~~~~~~~~~~~~~~~~~~~~~~~~~~~~~~~~~~~~~~~~~~~
(5.21a)\\[2mm] \end{eqnarray*}
where
$$l_{ij}=m_{ij}-(i-2\delta_{i,3}),\eqno(5.21b)$$ while $(m){_k}^{\pm ij}$
is a vector
obtained from $(m)$ by replacing $m_{ij}$ with $m_{ij}\pm 1$, they really
satisfy the
commutation relations (5.1a)--(5.1e) for $U_{q}[gl(2/1)_{0}]$. The highest
weight
described by the first row (signature)
$$[m]_{k}=[m_{12},m_{22},m_{32}]\eqno(5.22)$$ of the patterns (5.20) is nothing
but an ordered set of the eigen--values of the Cartan generators $E_{ii}$,
$i=1,2,3$,
on the highest weight vector $(M)_{k}$ defined as follows
$$E_{12}(M)_{k}=0,\eqno(5.23)$$
$$E_{ii}(M)_{k}=m_{i2}(M)_{k}, \eqno(5.24)$$ The highest weight vector
$(M)_{k}$
can be obtained from $(m)_{k}$ by setting $m_{11}=m_{12}$
$$(M)_{k}=\left[
\begin{array}{lcr}

\begin{array}{c}
m_{12}~~~m_{22}\\ m_{12}
\end{array}
;
\begin{array}{c}
m_{32}=m_{31}\\ m_{31}
\end{array}
\end{array}
\right].\eqno(5.25)$$
A lower weight vector $(m)_{k}$ can be derived {\it vice versa} from
$(M)_{k}$ by the
formula
$$(m)_{k}=\left({[m_{11}-m_{22}]!\over
[m_{12}-m_{22}]![m_{12}-m_{11}]!}\right)^{1/2} (E_{21})^{m_{12}-
m_{11}}(M)_{k}.\eqno(5.26)$$ Especially, for the case $k=0$, instead of the
above
notations, we skip the subscript 0, i.e.,
$$(m)_{0}\equiv (m);~~ [m]_{0}\equiv [m];~~ (M)_{0}\equiv (M), \eqno(5.27)$$
putting
$$m_{i2}=m_{i3},~~~ i=1,2,3,\eqno(5.28)$$ where $m_{i3}$ are some of the
complex values of $m_{i2}$, therefore, $m_{13}-m_{11},~ m_{11}-m_{23}\in {\bf
Z_{+}}$. We emphasize that $[m]$ and $(M)$, because of (5.7), are also,
respectively, the highest weight and the highest weight vector in the
$U_{q}[gl(2/1)]$--module $W^{q}=W^{q}([m])$. Characterizing the latter module as
the whole, $[m]$ and $(M)$ are, respectively, referred to as the global highest
weight and the global highest weight vector, while $[m]_{k}$ and $(M)_{k}$ are,
respectively, the local highest weights and the local highest weight vectors
characterizing only the submodules $V_{k}^{q}=V_{k}^{q}([m]_{k})$.\\

Following the arguments of Ref. \cite{k4}, for an alternative with (5.17)
basis of
$W^{q}$, we can choose the union of all the bases (5.20) which are denoted now
by the patterns $$\left[
\begin{array}{lcc}
m_{13}& m_{23}& m_{33} \\
m_{12}& m_{22}& m_{32}\\
m_{11}& 0 & m_{31}
\end{array}
\right]_{k}
\equiv
\left[
\begin{array}{lcr}

\begin{array}{c}
m_{12}~~~m_{22}\\ m_{11}
\end{array}
;
\begin{array}{c}
m_{32}=m_{31}\\ m_{31}
\end{array}
\end{array}
\right]_{k}\equiv (m)_{k}
,\eqno(5.29)$$
where the first row $[m]=[m_{13},m_{23},m_{33}]$ is simultaneously the highest
weight of the submodule $V_{0}^{q}=V_{0}^{q}([m])$ and the whole module
$W^{q}=W^{q}([m])$, while the second row $[m]_{k}=[m_{12},m_{22},m_{32}]$
is the local highest weight of some $U_{q}gl[(2/1)_{0}]$--module
$V^{q}_{k}=V^{q}_{k}([m]_{k})$ containing the considered vector $(m)_{k}$. The
basis (5.29) of $W^{q}$ is called the $U_{q}[gl(2/1)]$--reduced basis or
simply the
reduced basis. The latter representing a modified GZ basis description is
convenient for us
in investigating the module structure of $W^{q}$. Note once again that the
condition
$$m_{32}=m_{31}\eqno(5.20c)$$ has to be kept always.\\

The highest weight vectors $(M)_{k}$, now, in notation (5.29) have the form

$$(M)_{k}=\left[
\begin{array}{lcc}
m_{13}& m_{23}& m_{33} \\
m_{12}& m_{22}& m_{32}\\
m_{12}& 0 & m_{31}
\end{array}
\right]_{k},
\eqno(5.30)$$
as for $k=0$ the notations (5.27) and (5.28) are also taken into
account.\\[4mm]
{\bf Lemma 5.4}: {\it The highest weight vectors $(M)_{k}$ are expressed in
terms of
the induced basis (5.17) as follows} \begin{eqnarray*}
~~~~~~~~~~~~~~~~~~~~~~~~~~ (M)_{0}& = &a_{0}\left|0,0;(M)\right>,
~~~~a_{0}\equiv 1, \\[2mm]
(M)_{1}& = &a_{1}\left|0,1;(M)\right>, \\[2mm] (M)_{2}& =
&a_{2}\left\{\left|1,0,;(M)\right> +q^{2l}[2l]^{-1/2}
\left|0,1;(M)^{-11}\right>\right\},\\[2mm] (M)_{3}& =
&a_{3}\left\{\left|1,1;(M)\right>\right\},
~~~~~~~~~~~~~~~~~~~~~~~~~~~~~~~~~~~~~~~~~~~~~(5.31a) \end{eqnarray*}
{\it where $a_{i}$, $i=0,1,2,3$, are some numbers depending, in general, on
$q$,
while $l$ is}
$$l={1\over 2}(m_{13}-m_{23}).\eqno(5.31b)$$ \vspace*{2mm}
{\it Proof}: Indeed, all the vectors $(M)_{k}$ given above satisfy the
condition (3.4).\\


  From
formulae (5.24) and (5.31) the highest weights $[m]_{k}$ can be easily
identified
\begin{tabbing} \=123456791234567891234567890\= $[m]_{kk}$ \= =x \= $[m_{13}-
1,m_{23}-1,m_{33}+1,m_{43}+1]$,\=\kill \>\>
$[m]_{0}$ \> = \> $[m_{13}, m_{23}, m_{33}]$,\\[2mm] \>\>$[m]_{1}$ \> = \>
$[m_{13},
m_{23}-1, m_{33}+1]$,\\[2mm] \>\>$[m]_{2}$ \> = \> $[m_{13}-1, m_{23},
m_{33}+1]$,\\[2mm]
\>\>$[m]_{3}$ \> = \> $[m_{13},
m_{23}, m_{33}+2]$~~~~~~~~~~~~~~~~~~~~~~~~~~~~~~(5.32) \end{tabbing}
\vspace*{2mm}
Using the rule (5.26) we obtain all the basis vectors $(m)_{k}$:
\begin{eqnarray*}
~~~~~~~~~~~~~~~~~~~~~~~~~~~~ (m)_{0}& \equiv & \left[
\begin{array}{lcc}
m_{13}& m_{23}& m_{33} \\
m_{13}& m_{23}& m_{33}\\
m_{11}& 0 & m_{33}
\end{array}
\right]
=\left|0,0,;(m)\right>,\\[4mm]
(m)_{1}&\equiv &
\left[
\begin{array}{lcc}
m_{13}& m_{23}& m_{33} \\
m_{13}& m_{23}-1& m_{33}+1\\
m_{11}& 0 & m_{33}+1
\end{array}
\right] \\[4mm]
& = &a_{1}\left\{-\left(
{[l_{13}-l_{11}]\over
[2l+1]}\right)^{1/2}\left|1,0;(m)^{+11}\right>\right. \\ & &\left.
+q^{2(l_{11}-l_{13})}\left(
{[l_{11}-l_{23}]\over
[2l+1]}\right)^{1/2}\left|0,1;(m)\right>\right\}, \\[4mm]
(m)_{2}&\equiv &
\left[
\begin{array}{ccc}
m_{13}& m_{23}& m_{33} \\
m_{13}-1& m_{23}& m_{33}+1\\
m_{11}& 0 & m_{33}+1
\end{array}
\right] \\[4mm]
& = &a_{2}\left\{\left(
{[l_{11}-l_{23}]\over
[2l]}\right)^{1/2}\left|1,0;(m)^{+11}\right>\right. \\ & &\left.
+q^{l_{11}-l_{23}}\left(
{[l_{13}-l_{11}]\over
[2l]}\right)^{1/2}\left|0,1;(m)\right>\right\}, \\[4mm]
(m)_{3}&\equiv &
\left[
\begin{array}{ccc}
m_{13}& m_{23}& m_{33} \\
m_{13}-1& m_{23}-1& m_{33}+2\\
m_{11}& 0 & m_{33}+2
\end{array}
\right] \\[4mm]
&=& a_{3}\left|1,1;(m)\right>,
~~~~~~~~~~~~~~~~~~~~~~~~~~~~~~~~~~~~~~~~~~~~~~~
(5.33) \end{eqnarray*}
where $l_{ij}$ and $l$ are given in (3.21b) and (5.31b), respectively.
Here, we skip
the subscript $k$ in the patterns given above since there are no degenerations
between them. The formulae (5.33), in fact, represent the way in which the
reduced
basis (5.29) is written in terms of the induced basis (5.16). From (5.33)
we can
derive their invert relation \begin{eqnarray*}
~~~~~~~~~~~~~~~~~~~~~~~~~~~~~~ \left|0,0;(m)\right>& = &(m)_{0}\equiv
(m)\\[2mm]
\left|1,0;(m)\right>& = &-{1\over a_{1}}q^{l_{11}-l_{23}-1}
\left({[l_{13}-l_{11}+1]\over
[2l+1]}
\right)^{1/2}(m)_{1}^{-11}\\[2mm] & & +{1\over a_{2}}q^{l_{11}-l_{13}-1}
{\left([l_{11}-
l_{23}-1][2l] \right)^{1/2} \over [2l+1]}(m)_{2}^{-11},\\[2mm]
\left|0,1;(m)\right>& = &{1\over a_{1}}
\left({[l_{11}-l_{23}]\over [2l+1]}
\right)^{1/2}(m)_{1}\\[2mm] & & +{1\over a_{2}}
{\left([l_{13}-l_{11}][2l]\right)^{1/2}
\over [2l+1]} (m)_{2},\\[2mm]
\left|1,1;(m)\right>& = &{1\over
c_{3}}(m)^{-11}_{3}.
~~~~~~~~~~~~~~~~~~~~~~~~~~~~~~~~~~~~~~~~~~~~~(5.34) \end{eqnarray*}

Now we are ready to compute all the matrix elements of the generators in
the basis
(5.29). As we shall see, the latter basis allows an evident description of a
decomposition of a $U_{q}[gl(2/1)]$--module
$W^{q}$ in irreducible $U_{q}[gl(2/1)_{0}]$--modules $V^{q}_{k}$. Since the
finite--dimensional representations of the $U_{q}[gl(2/1)]$ in some basis
are completely
defined by the actions of the even generators and the odd Weyl--Chevalley ones
$E_{23}$ and $E_{32}$ in the same basis, it is sufficient to write down the
matrix
elements of these generators only. For the even generators the matrix elements
have already been given in (5.21), while for $E_{23}$ and $E_{32}$, using the
relations (5.1)--(5.3), (5.33) and (5.34) we have
\begin{eqnarray*}
~~~~~~~~~~~~~~~~~~~~~~E_{23}(m) &=& 0,\\[2mm]
E_{23}(m)_{1} &=&a_{1}
\left({[l_{11}-l_{23}]\over [2l+1]}\right)^{1/2}[l_{23}+l_{33}+3] (m),\\[4mm]
E_{23}(m)_{2} &=&a_{2}
\left({[l_{13}-l_{11}]\over [2l]}\right)^{1/2}[l_{13}+l_{33}+3] (m),\\[4mm]
E_{23}(m)_{3} &=&a_{3}
\left\{{1\over a_{1}q}\left({[l_{13}-l_{11}]\over
[2l+1]}\right)^{1/2}[l_{13}+l_{33}+3](m)_{1}\right. \\[2mm] &&\left.
-{1\over a_{2}q}\left({[l_{11}-l_{23}]
[2l]}\right)^{1/2}{[l_{23}+l_{33}+3]\over [2l+1]}(m)_{2}\right\}
~~~~~~~~~~~~~(5.35a)
\end{eqnarray*}
and
\begin{eqnarray*}
~~~~~~~~~~~~~~~~~~~~~~~~~~~~~~~~~~~~~ E_{32}(m) &=& {1\over
a_{1}}\left({[l_{11}-
l_{23}]\over [2l+1]}\right)^{1/2}(m)_{1}\\[2mm]
&&+{1\over a_{2}}
{\left([l_{13}-l_{11}][2l]\right)^{1/2}\over [2l+1]}(m)_{2}\\[4mm]
E_{32}(m)_{1} &=&{a_{1}\over
a_{3}}q
\left({[l_{13}-l_{11}]\over [2l+1]}\right)^{1/2} (m)_{3},\\[4mm]
E_{32}(m)_{2} &=&-{a_{2}\over a_{3}}q
\left({[l_{11}-l_{23}]\over [2l]}\right)^{1/2} (m)_{3},\\[4mm]
E_{32}(m)_{3} &=&0.
~~~~~~~~~~~~~~~~~~~~~~~~~~~~~~~~~~~~~~~~~~~~~~~~~ (5.35b)
\end{eqnarray*}
{\bf Lemma 5.5}: {\it The finite--dimensional representations (5.35) of
$U_{q}[gl(2/1)]$ are irreducible and called typical if and only if the
condition}
$$[l_{13}+l_{33}+3][l_{23}+l_{33}+3]\neq 0\eqno(5.36)$$ {\it holds}.\\
{\it Proof}: By the
same argument useed in Ref. {k4} we can conclude that $W^{q}$ is
irreducible if and only if
$$E_{14}E_{23}E_{13}E_{31}E_{32}E_{41}\otimes (M)
\neq 0.$$
The latest condition in turn can be
proved, after some elementary calculations, to be equivalent to
$$[E_{11}+E_{33}+1][E_{22}+E_{33}](M)
\neq 0,$$
\vspace*{1mm}
which is nothing but the condition (5.36).\\

In case the condition (5.36) is violated, i.e. one of the following
condition pairs
$$[l_{13}+l_{33}+3]=0\eqno(5.37a)$$
and
$$[l_{23}+l_{33}+3]\neq 0\eqno(5.37b)$$
or
$$[l_{13}+l_{33}+3]\neq 0\eqno(5.38a)$$
and
$$[l_{23}+l_{33}+3]=0\eqno(5.38b)$$
(but not both (5.37a) and (5.38b) simultaneously) holds, the module $W^{q}$
is no longer irreducible but indecomposable. However, there exists an
invariant subspace, say $I_{k}^{q}$, of $W^{q}$ such that the
factor--representation
in the factor--module $$W_{k}^{q}:=W^{q}/I_{k}^{q}\eqno(5.39)$$ is
irreducible. We say that is a nontypical representation in a nontypical
module $W_{k}^{q}$.
Then, as in
Ref. \cite{k5}, it is not difficult for us to prove the following
assertions\\[4mm]
{\bf Lemma 5.6}:
$$V_{3}^{q}\subset I_{k}^{q},\eqno(5.40)$$ {\it and}
$$V_{0}^{q}\cap I_{k}^{q}=\emptyset .\eqno(5.41)$$
>From (5.35)--(5.38) we can
easily find all nontypical representations of $U_{q}[gl(2/1)]$ which are
classified in 2
classes.\\[5mm]
\noindent
{\bf 5.b. Nontypical representations of U$_{q}$[gl(2/1)]}\\

1) \underline{Class 1 nontypical representations}:\\

This class is characterized by the conditions $(5.37a)$ and $(5.37b)$ which for
generic $q$ take the forms $$l_{13}+l_{33}+3=0,\eqno(5.37x)$$
and
$$l_{23}+l_{33}+3\neq 0,\eqno(5.37y)$$
respectively.
In other words, we have to replace everywhere all $m_{33}$ by $-m_{13}-1$
and keep $(5.37y)$ valid. Thus we have\\[4mm]
{\bf Lemma 5.7}: {\it
The class 1
maximal invariant subspace in} $W^{q}$ {\it is}
$$I_{1}^{q}=V_{3}^{q}\oplus V_{2}^{q}.\eqno(5.42)$$
{\it Proof}: Applying (5.37) to (5.35) we obtain (5.42).\\

 Then the class 1 nontypical representations in
$$W_{1}^{q}=W_{1}^{q}([m_{13},m_{23},-m_{13}-1])\eqno(5.43)$$ are given
through (5.35) by keeping the conditions (5.37) (i.e., $(5.37x)$ and
$(5.37y)$) and
replacing with 0 all vectors belonging to $I_{1}^{q}$: \begin{eqnarray*}
~~~~~~~~~~~~~~~~~~~~~~~~~~~~~ E_{23}(m)&=&0,\\[2mm]
E_{23}(m)_{1}&=&a_{1}
\left({[l_{11}-l_{23}]\over [2l+1]}\right)^{1/2}[l_{23}-l_{13}](m)
~~~~~~~~~~~~~~~~~~~~~~(5.44a)
\end{eqnarray*}
and
\begin{eqnarray*}
~~~~~~~~~~~~~~~~~~~~~~~~~~~~~~~~~~~~
E_{32}(m)&=&{1\over a_{1}}\left({[l_{11}-l_{23}] \over
[2l+1]}\right)^{1/2}(m)_{1},\\[2mm]
E_{32}(m)_{1}&=&0.~~~~~~~~~~~~~~~~~~~~~~~~~~~~~ ~~~~~~~~~~~~~~~~~~~~~~~~(5.44b)
\end{eqnarray*}

2) \underline{Class 2 nontypical representations}:\\

For this class nontypical representations we must keep the conditions

$$l_{13}+l_{33}+3\neq 0,\eqno(5.38x)$$
and
$$l_{23}+l_{33}+3= 0.\eqno(5.38y)$$
derived, respectively, from $(5.38a)$ and $(5.38b)$ when the deformation
parameters $q$ are generic.
Equivalently, we have to replace everywhere all $m_{33}$ by $-m_{23}$ and keep
$(5.38x)$ valid.\\

Now the invariant subspace $I_{2}^{q}$ is the following \\[4mm]
{\bf Lemma 5.8}:
{\it The class 2 maximal invariant subspace in} $W^{q}$ {\it is}
$$I_{2}^{q}=V_{3}^{q}\oplus V_{1}^{q}.\eqno(5.45)$$
{\it Proof}: Using (5.38) in (5.35) we derive (5.45).\\

 The class 2 nontypical representations in
$$W_{2}^{q}=W_{2}^{q}([m_{13},m_{23},- m_{23}])\eqno(5.46)$$ are also
given through (5.35) but by keeping the conditions
(5.38) (i.e., $(5.38x)$ and $(5.38y)$) valid and replacing by 0 all
vectors belonging to the invariant subspace
$I_{2}^{q}$:
\begin{eqnarray*}
~~~~~~~~~~~~~~~~~~~~~~~~~~~~ E_{23}(m)&=&0,\\[2mm]
E_{23}(m)_{2}&=&a_{1}
\left({[l_{13}-l_{11}]\over [2l]}\right)^{1/2}[2l+1](m)
~~~~~~~~~~~~~~~~~~~~~~~~(5.47a)
\end{eqnarray*}
and
\begin{eqnarray*}
~~~~~~~~~~~~~~~~~~~~~~~~~~~~~~~~~~~~~~ E_{32}(m)&=&{1\over
a_{2}}{\left([l_{13}-l_{11}][2l]\right)^{1/2} \over [2l+1]}(m)_{2},\\[2mm]
E_{32}(m)_{2}&=&0.
~~~~~~~~~~~~~~~~~~~~~~~~~~~~~~~~~~~~~~~~~~~~~~~~~~ (5.47b)
\end{eqnarray*}

We have just considered the quantum superalgebra $U_{q}[gl(2/1)]$ and
constructed
all its typical and nontypical representations leaving the coefficients
$a_{i}$,
$i=1,2,3$, as free parameters which can be fixed by some additional
conditions, for
example, the Hermiticity condition. As an intermediate step (which,
however, is of
an independent interest) we also introduced the reduced basis (5.29) which,
as it is
an extension of the Gel'fand--Zetlin basis to the present case, is
appropriate for an
evident description of decompositions of $U_{q}[gl(2/1)]$--modules in
irreducible
$U_{q}[gl(2/1)_{0}]$--modules. We can prove the following propositions\\[4mm]
{\bf Lemma 5.9}: {\it The class of the finite--dimensional representations
determined
in this paper (Subsects. 5.a and 5.b), contains all
finite--dimensional irreducible representations of $U_{q}[gl(2/1)]$ and
$U_{q}[sl(2/1)]$}.\\[4mm]
and\\[4mm]
{\bf Lemma 5.10}: {\it The finite--dimensional representations of the quantum
superalgebra $U_{q}[gl(2/1)]$ are quantum deformations of the
finite--dimensional representations of the superalgebra $gl(2/1)$.}\\

The {\bf Lemma 5.9} is proved by similar arguments
as those used in the proofs of Proposition 9 and Proposition 10 in Ref.
\cite{k5}, while {\bf Lemma 5.10} can be verified by direct computations.\\

  Since the nontypical representations have only been well investigated
for a few cases of both classical and quantum superalgebras, the present
results can be considered as a small step forward in this
direction.\\[9mm] {\Large {\bf 6.
Conclusion}}\\

Certainly, many questions remain unconsidered in the framework of the present
paper but we hope that the latter gives relevant information about the
classical and
quantum superalgebras as well as an idea on their representations. Based on
Kac's representation theory for the superalgebras \cite{kac,kac3,kac4,sch}
we succeeded in finding all finite--dimensional representations and a wide
class of infinite--dimensional representations of several higher rank
superalgebras of nonparticular
types \cite{k3,k1,k2,neli}. Recently, extending that classical theory to the
quantum deformation we worked out a method for explicit constructions of
representations of quantum superalgebras \cite{k4,k6,k7,k8,k5}. Our method,
avoiding the use of the Clebsch--Gordan
coefficients which are usually unknown for higher rank (classical and quantum)
algebras, is applicable not only to the one--parametric deformations but
also to the multi--parametric ones (see, for example, Refs.
\cite{k6,k7,k8}). Moreover, our
approach may have an advantage
as it is worthy to mention that the theory of representations and,
especially, of the
nontypical representations is far from being complete even for the nondeformed
superalgebras. In particular, the dimensions of the nontypical
representations are
unknown unless the ones for $sl(1/n)$ computed recently in Ref.
\cite{schlo}. Based
on the generalizations of the concept of the GZ basis (see Refs.
\cite{k4},\cite{k5} and \cite{palev3} and references therein) the
matrix elements of all nontypical representations were computed only for
superalgebras of lower ranks or of particular types like $sl(1/n)$ and
$gl(1/n)$ in
Ref. \cite{palev4,palev6}. Later, the essentially typical representations of
$gl(m/n)$ were
also constructed \cite{palev5}. So far, however, the GZ basis concept was not
defined and presumably cannot be defined for nontypical $gl(m/n)$--modules with
$m,n \geq 2$. This was of why to try to describe the nontypical modules in
terms of
the basis of the even subalgebras. This approach, developed so far for
classical
superalgebras \cite{k3,k1,k2} and for quantum superalgebras \cite{k4,k6,k5})
turns out to be appropriate for explicit descriptions of all nontypical
modules of
$gl(2/2)$ (see Refs. \cite{k1,neli}), $U_{q}[gl(2/2)]$ (see Refs. \cite{k4} and
\cite{k5}) and multi--parametric quantum superalgebras \cite{k6,k7,k8}.
Our approach in Refs
\cite{k4,k6,k7,k8,k5}, unlike some earlier approaches, avoids, however,
the use of the
Clebsch--Gordan coefficients which are not always known for higher rank
(quantum
and classical) algebras. Other extensions were made in Ref. \cite{palev2}
for all
finite-dimensional representations of $U_{q}[gl(1/n)]$ and in Ref.
\cite{palev1} for a
class of finite--dimensional representations of $U_{q}[gl(m/n)]$. To the
best of our
knowledge, we gave for the first time \cite{k3,k1,k4,k6,k2,k5,neli}, explicit
expressions for all finite--dimensional representations or a wide class of
infinite--dimensional representations of several classical and quantum
superalgebras including those of higher ranks.\\

  We hope that our approach allowing to establish in consistent ways
defining relations of quantum (super)algebras \cite{k4,k7,k8} can be extended
to the case of deformation parameters being roots of unity (for
representations of quantum groups at roots of unity, see, for example,
\cite{kac2}).\\[9mm]
{\Large {\bf Acknowledgements}}\\

  I am grateful to Professor J. Tran Thanh Van and {\it "Rencontres du
Vietnam"} for financial support and Professor S. Randjbar--Daemi for kind
hospitality at the Abdus Salam International Centre for Theoretical
Physics, Trieste, Italy. I would like to thank the organizers of the {\it
V--th National Mathematics Conference} (Hanoi, 17--20 September 1997)
for giving me the opportunity to report on the present topic.\\

  This work was also suported by the Vietnam National Basic Research
Programme in Natural Science under the grant KT 4.1.5.\\


\end{document}